    \input amstex
     \documentstyle{amsppt}
   \pageno=1
\NoRunningHeads

     \magnification = 1200
    \pagewidth{5.3 true in}
    \pageheight{8.5 true in}
   \hoffset .6 true in

   \NoBlackBoxes
     \parskip=3 pt
     \parindent = 0.3 true in

\def\({\bigl(}
\def\){\bigr)}
\def\Z{\Bbb Z}

\def\Q{\Bbb Q}
\def\R{\Bbb R}
\def\C{\Bbb C}

\def\G{\Cal G}
\def\R{\Cal R}
\def\S{\Cal S}

\def\N{\operatorname{N}}

\def\Oc{\Cal O}

\def\fa{\frak a}

\def\ep{\varepsilon}
\def\p{\frak p}
\def\m{\frak m}

\def\m{\frak m}
\def\epl{\frac{1+\tau}{2}}
\def\emi{\frac{1-\tau}{2}}
\def\Fit{{\operatornamewithlimits{Fit}}}
\def\Ann{{\operatornamewithlimits{Ann}}}

\def\FI{{\Cal I}^{\text{Fit}}}
\def\SI{{\Cal I}^{\text{Sti}}}

\def\S{{\Cal S}}

\def\OE{\Oc_E^S}
\def\OF{\Oc_F^S}
\def\KE{K_2(\OE)}
\def\KF{K_2(\OF)}
\def\WE{W_2(E)}
\def\WF{W_2(F)}
\def\AW{1}
\def\JTT{13}
\def\Ga{3}
\def\QuFG{10}
\def\MW{9}
\def\Wi{14}
\def\Ko{7}
\def\KoMo{8}
\def\Ka{6}
\def\SaSi{11}
\def\Gr{4}
\def\Gre{5}
\def\Sch{12}
\def\Eis{2}
    \topmatter
    \title
Values at $s=-1$ of $L$-functions \\
for relative quadratic extensions \\
of number fields, and \\
the Fitting ideal of the tame kernel
    \endtitle
\author
Jonathan W. Sands
\\
Dept. of Mathematics and Statistics
\\
University of Vermont, Burlington VT 05401 USA
\\
Jonathan.Sands\@uvm.edu
\endauthor

\thanks
2000 {\it MSC:} Primary 11R42; Secondary 11R70, 19F27.
\newline\indent Thanks to the Department of Mathematics at the University of California, San
\newline\indent Diego for its hospitality during a 2006--2007 sabbatical visit when this paper
\newline \indent was written.
\endthanks

\abstract
  Fix a relative quadratic extension $E/F$ of totally real number fields and let $G$ denote the
Galois group of order 2. Let $S$ be a finite set of primes of $F$
containing the infinite primes and all those which ramify in $E$,
let $S_E$ denote the primes of $E$ lying above those in $S$, and let
$\Oc_E^S$ denote the ring of $S_E$-integers of $E$. Assume the truth
of the 2-part of the Birch-Tate conjecture relating the order of the
tame kernel $K_2(\Oc_E)$ to the value of the Dedekind zeta function
of $E$ at $s=-1$, and assume the same for $F$ as well. We then prove
that the Fitting ideal of $K_2(\Oc_{E}^S)$ as a $\Z[G]$-module is
equal to a generalized Stickelberger ideal. Equality after tensoring
with $\Z[1/2][G]$ holds unconditionally.

\endabstract

    \endtopmatter

     \document
\baselineskip=20pt

\bigskip\bigskip\bigskip
 \vfill\eject

\heading{I. Introduction}
\endheading

Let $E/F$ be a fixed relative quadratic extension of algebraic
number fields with Galois group $G=\langle \tau \rangle$. We also
fix a finite set $S$ of primes of $F$ which contains all of the
infinite primes of $F$ and all of the primes which ramify in $E$.
Associated with this data is a Stickelberger function, or
equivariant $L$-function, $\theta_{E/F}^S(s)$. It is a meromorphic
function of $s$ with values in the group ring $\C[G]$. To define it,
let $\p$ run through the (finite) primes of $F$ not in $S$, and
$\fa$ run through integral ideals of $F$ which are relatively prime
to each of the elements of $S$. Then ${\N} \fa$ denotes the absolute
norm of the ideal $\fa$,
 $\sigma_{\fa} \in G$ is the well-defined automorphism
attached to $\fa$ via the Artin map, and

$${
\theta_{E/F}^S(s)=
\sum\Sb {\frak a} \ \text{integral} \\ ({\frak a},S)=1 \endSb
\dfrac{1}{ {\N} {\frak a}^s}\sigma_{\fa}^{-1}=
\prod\Sb \text{prime } \frak p \notin S \endSb
\bigl(1-\dfrac{1}{ {\N} {\frak p}^{s}}} \sigma_{\frak p}^{-1}\bigr)^{-1}.$$

 These expressions converge for the real part of $s$ greater than 1 and the function they define
extends meromorphically to all of $\C$.

The function $\theta_{E/F}^S(s)$ is connected with the arithmetic of
the number fields $E$ and $F$ in ways one would like to make as
precise as possible. Define the ring of $S$-integers $\OF$ of $F$ to
be the set of elements of $F$ whose valuation is non-negative at
every prime not in $S$. Similarly, define the ring $\OE$ of
$S$-integers of $E$ to be the set of elements of $E$ whose valuation
is non-negative at every prime not in $S_E$, the set of all primes
of $E$ which lie above some prime in $S$. The zeta-functions defined
by $\zeta_F^S(s) \sigma_{\Oc_F}=\theta_{F/F}^S(s)$ and $\zeta_E^S(s)
\sigma_{\Oc_E}=\theta_{E/E}^{S_E}(s)$ may be viewed as the
zeta-functions of the Dedekind domains $\OF$ and $\OE$.

Our focus will be on the ``higher Stickelberger element" $\theta_{E/F}^S(-1)$.
It is conjecturally related to the algebraic $K$-group $\KE$. This
group is known to be finite by [\Ga] and [\QuFG], and called the
tame kernel of $E$.

  The precise statement of the conjectured arithmetic interpretation
of $\theta_{E/F}^S(-1)$ will also involve another finite group.
Let $\mu_\infty$ denote the group of all roots of unity in an algebraic closure
$\overline{\Q}$ of $\Q$ containing $E$, and let $\G$ denote the
Galois group of $\overline{\Q}/\Q$. Define $W_2=W_2(\overline{\Q})$
to be the $\Z[\G]$-module whose underlying group is $\mu_\infty$,
with the action of $\gamma \in \G$ on $\omega \in W_2$ given by
$\omega^\gamma=\gamma^2(\omega)$. Then for any subfield $L$ of
$\overline{\Q}$, let $W_2(L)$ be the submodule fixed under this
action by the Galois group of $\overline{\Q}$ over $L$. Then
$W_2(E)$ naturally becomes a $\Z[G]$-module, where the action of $G$
arises by lifting elements of $G$ to the Galois group of
$\overline{\Q}$ over $F$. One easily sees that the $G$-fixed
submodule $W_2(E)^G$ equals $W_2(F)$. We use the notation
$w_2(L)=|W_2(L)|$, which we note is finite for any algebraic
number field $L$.

We are now ready to state the conjecture of Birch and Tate (see
section 4 of [\JTT]), which gives a precise arithmetic
interpretation of $\zeta_F^S(-1)$. We state an extended form of it
for arbitary finite $S$ which is an easy consequence of the original
conjecture for minimal $S$ (see Corollary 3.3 of [\SaSi]).

\proclaim{Conjecture 1.1 (Extended Birch-Tate)}
Suppose that $F$ is totally real. Then
$$\zeta_F^S(-1)= (-1)^{|S|} \dfrac{|\KE|}{w_2(F)}$$
\endproclaim

 Deep results on Iwasawa's Main conjecture in [\MW] and [\Wi] lead to the following
(see [\Ko]).

\proclaim{Theorem 1.2}
The Birch-Tate Conjecture holds if $F$ is abelian over $\Q$, and the odd part
holds for all totally real $F$.
\endproclaim

We note that the 2-part of the Birch-Tate conjecture for $F$ would
follow from the 2-part of Iwasawa's Main conjecture for $F$.
We can now state our main result,
using the annihilator $\Ann_{\Z[G]}(W_2(E))$ of $W_2(E)$ in
$\Z[G]$.

\proclaim{Theorem 1.3} Let $E/F$ be a relative quadratic extension of
totally real number fields, with Galois group $G$.
Assume that the 2-part of the Birch-Tate conjecture holds for $E$ and for $F$.
Then the (first) Fitting ideal of $\KE$ as a $\Z[G]$-module is
$$\Fit_{\Z[G]}(\KE)=\Ann_{\Z[G]}(W_2(E))\theta_{E/F}^S(-1).$$
Equality of the extension of these ideals to $\Z[1/2][G]$ holds unconditionally.
\endproclaim

  This result is to be compared with results computing the Fitting
ideals of certain modified ideal class groups in terms of more
classical Stickelberger ideals. Notable among these is the recent
result of Greither [\Gre] identifying the Fitting ideal of the dual
of the odd part of the minus part of the ideal class group of a CM
field which is abelian over a totally real field. Assuming the
Equivariant Tamagawa Number Conjecture, this ideal is shown to equal
a Stickelberger ideal obtained from the values at $s=0$ of
Stickelberger functions for subextensions. Our Theorem 1.3 begins to
provide an analog at $s=-1$. The comparison can provide some insight
since the odd part in Theorem 1.3 is unconditional, and the result
on the 2-part, though conditional, has no counterpart as yet at
$s=0$.

For ease of notation from now on, let $\R=\Z[G]$. Denote the
Fitting ideal we
are interested in by $\FI=\Fit_{\R}(\KE)$, and the generalized
Stickelberger ideal by $\SI=\Ann_{\R}(W_2(E))\theta_{E/F}^S(-1)$.
 The ingredients we will need to prove
their equality are cohomology computations, properties of Fitting
ideals, and some information on the ideal structure of $\R$ and
related rings.

\heading{II. Cohomology of $K_2(\Oc_E^S)$}
\endheading

 We continue to assume throughout that $E/F$ is a relative quadratic extension of
totally real fields with Galois group $G=\langle \tau \rangle$. For
a $\Z[G]$-module $M$, we define $M^G$ as usual to be the submodule of
$M$ which is annihilated by $1-\tau$, and $M^-$ to be the submodule
annihilated by $1+\tau$. Also let $M_G$ denote the module of
co-invariants of $M$ under the action of $G$. The following result
provides the key to computing the orders we will need.

\proclaim{Theorem 2.1} Under our assumptions, the transfer map from $\KE$ to $\KF$ is surjective
with kernel $K_2(\Oc_E^S)^{1-\tau}$. So we have a short exact sequence
$$0 \rightarrow K_2(\Oc_E^S)^{1-\tau}\rightarrow \KE \overset{\text{Trans}} \to \rightarrow \KF \rightarrow 0 $$
\endproclaim

\demo{Proof} In the more general setting of a Galois extension of number fields $E/F$, with Galois group $G$,
 Kahn's theorem 5.1 of [\Ka] leads to an exact sequence induced by the transfer map and involving the number
$r_\infty(E/F)$  of infinite primes of $F$ which ramify in $E$ (see
[\KoMo, Prop. 1.6] or [\SaSi, Thm. 3.4]):
$$0 \rightarrow K_2(\Oc_E^S)_G\rightarrow K_2(\Oc_F^S)\rightarrow \{\pm 1\}^{r_\infty(E/F)}\rightarrow 0 $$
In our situation with $G=\langle \tau \rangle$ cyclic and $E$ totally real, $K_2(\Oc_E^S)_G=\KE/\KE^{1-\tau}$ and
$r_\infty(E/F)=0$. The result follows.
\qed
\enddemo

\proclaim{Proposition 2.2} $|\KE^G|=|\KF|$
\endproclaim
\demo{Proof} Theorem 2.1 gives $|\KF|=|\KE|/|\KE^{1-\tau}|$, while
the exact sequence
$$0 \rightarrow K_2(\Oc_E^S)^G\rightarrow \KE \overset{1-\tau} \to \rightarrow \KE^{1-\tau} \rightarrow 0 $$
gives the same value for $|\KE^G|$. \qed
\enddemo

Let
 $\Q_r$ denote the $r$th layer of the cyclotomic
$\Z_2$ extension of $\Q$; it may be defined as
$\Q_r=\Q(\mu_{2^{r+2}})^+$, the maximal real subfield of the field
of $2^{r+2}$th roots of unity. Thus $\Q_0=\Q$ and
$\Q_1=\Q(\sqrt{2})$. Also let $\Q_\infty$ denote the union of the
$\Q_r$ for all natural numbers $r$. This is the cyclotomic
$\Z_2$-extension of $\Q$.
The cyclotomic $\Z_2$-extension of a field $L$ is
$L_\infty=L\cdot \Q_\infty$, and the $n$-th layer of this extension
$L_n$ is the unique subfield of $L_\infty$ of degree $2^n$ over $L$.

 \proclaim{Proposition 2.3} The cohomology groups
$$H^1(G,\KE)=\KE^-/\KE^{1-\tau}$$ and
$$H^2(G,\KE)=\KE^G/\KE^{1+\tau}$$ both have order 1 if $E=F_1$ is the first layer in the
cyclotomic $\Z_2$-extension of $F$, and both have order 2 otherwise.
\endproclaim

\demo{Proof} Since $\KE$ is finite and $G$ is cyclic, both
cohomology groups have the same order (see [\AW, Prop.11]). To
compute the order of \newline
$H^2(G,\KE)=\KE^G/\KE^{1+\tau}$, view it as the
cokernel of the map $\KE\overset{1+\tau}\to \rightarrow \KE^G$. A
standard functorial property states that this map factors as the
transfer map $\KE \overset{\text{Trans}} \to \rightarrow \KF$
followed by the map $\KF \overset{\iota_*} \to \rightarrow \KE^G $
induced by the inclusion $\OF \overset{\iota} \to \rightarrow \OE$
of rings. Since the transfer is surjective, by Theorem, 2.1, we are
reduced to computing the order of the cokernel of $\iota_*$. By
Proposition 2.2, the domain and codomain of this map have the same
order, and thus the order of the cokernel equals the order of the
kernel. The order of $\ker(\iota_*)$ is computed in [\SaSi, Lemma
7.3], yielding 1 if $E=F_1$ and 2 otherwise. \qed
\enddemo

\heading{III. Cohomology of $\WE$}
\endheading

 We begin with two simple lemmas.

\proclaim{Lemma 3.1} If $L$ is a real field, and $L\cap \Q_\infty =
\Q_r$, then $2^{r+3}$ exactly divides $w_2(L)$. In particular,
$w_2(L)$ is always divisible by 8.
\endproclaim
\demo{Proof} This follows from the fact that $L(\mu_{2^{r+2+k}})$ is
the composite of the cyclic extensions $L(\mu_4)$ and $L\cdot
\Q_{r+k}$. See [\SaSi, Lemma 7.2] for full details.
\qed
\enddemo

\proclaim{Lemma 3.2} Define the positive integers $r$ and $s$ by
$F\cap \Q_\infty =\Q_r$, and $E\cap \Q_\infty =\Q_s$. Then $s>r$ if
and only if $E=F_1$, the first layer of the cyclotomic
$\Z_2$-extension of $F$. (In fact, $s\leq r+1$ always holds.)
\endproclaim
\demo{Proof}
We have $s>r$ if and only if $E\supset \Q_{r+1}\not\subset F$
if and only if $E\supset F\cdot\Q_{r+1}=F_1\supset F$.
Since $[E:F]=2=[F_1:F]$, this last condition is
equivalent to $E=F_1$.
\qed
\enddemo

For any prime number $p$, let $\Z_{(p)}\subset \Q$ denote the localization
of $\Z$ at the prime ideal $(p)$. If $M$ is a $\Z$-module, then similarly
let $M_{(p)}\cong M\otimes_{\Z}\Z_{(p)}$ be the localization of $M$
at $(p)$.  Note that if $M$ is a finite abelian group, we may identify
$M_{(p)}$ with the Sylow $p$-subgroup of $M$,
which we call the $p$-part of $M$. If $M$ is a $\Z[G]$-module, then
$M_{(p)}$ becomes a $\Z_{(p)}[G]$-module. If $I$ is an ideal of $\R=\Z[G]$, we may identify
$I_{(p)}$ with the extended ideal $I\Z_{(p)}[G]\subset \Z_{(p)}[G]=\R_{(p)}\subset \Q[G]$.

\proclaim{Corollary 3.3} We have $\WE_{(2)} = \WF_{(2)}$ if and only
if $E \neq F_1$.
\endproclaim
\demo{Proof}
Since $\WE$ is cyclic and contains $\WF$, this follows upon applying the two preceeding lemmas.
\qed
\enddemo

\proclaim{Proposition 3.4} The cohomology groups $$H^1(G,\WE)=\WE^-/\WE^{1-\tau}$$ and
$$H^2(G,\WE)=\WE^G/\WE^{1+\tau}=\WF/\WE^{1+\tau}$$ have order 1 if $E=F_1$ is the first layer in the
cyclotomic $\Z_2$ extension of $F$, and order 2 otherwise.
\endproclaim
\demo{Proof}
Since $\WE$ is finite and $G$ is cyclic, the two cohomology groups have the same order.
We will compute $H^2(G,\WE)$.
Since $|G|=2$, the cohomology groups have exponent 2.
Expressing $\WE$ as a direct product (as a $\Z[G]$-module) of its 2-part $\WE_{(2)}$
and its odd part, we therefore need only consider the cohomology of $\WE_{(2)}$.

According to Lemma 3.1, we may write $w_2(F)=8m$. The action of
$\tau$ on the cyclic group $\WE$ is non-trivial and fixes $\WF$.
Hence $\tau$ must act as $1+kw_2(F)$ for some integer $k$. Thus
$\WE_{(2)}^{1+\tau}=\WE_{(2)}^{2+kw_2(F)}=(\WE_{(2)}^{1+4m})^2=\WE_{(2)}^2$.
Since $\WE^G=\WF$, we conclude that
$H^2(G,\WE_{(2)})=\WF_{(2)}/\WE_{(2)}^2$. Again using the cyclicity
of $\WE$, we see that this cohomology group is cyclic of exponent 2
and is non-trivial if and only if $\WE_{(2)}= \WF_{(2)}$. By the
corollary, this is equivalent to $E\neq F_1$. (Here we also see that
$s\leq r+1$.) Thus the cohomology groups are of order 2 when $E\neq
F_1$, and are trivial otherwise.
\qed
\enddemo

\heading{IV. } The Annihilator of $\WE$
\endheading

We have observed that
$\tau$ acts on $\WE$ as $N=1+kw_2(F)$, for some integer $k$. Thus
$\tau-N$ and $|\WE|$ lie in $\Ann_\R(\WE)$. It is easy to see that
these two elements generate $\Ann_\R(\WE)$. For if $a\tau+b\in
\Ann_\R(\WE)$, then $a\tau+b -a(\tau - N)=b+aN\in \Ann_\R(\WE)$. We
see that $(b+aN)=c|\WE|$, for some $c\in \Z$ and so
$a\tau+b=a(\tau-N)+c|\WE|$.

For future purposes, we would like to express $\Ann_\R(\WE)$
in terms of $1+\tau$ and $1-\tau$.

\proclaim{Proposition 4.1}
$\Ann_\R(\WE)$ is generated by $|\WE^{1+\tau}|(1+\tau)$ and
$|\WE^{1-\tau}|(1-\tau)$ when $E\neq F_1$. It is generated by
$|\WE^{1+\tau}|\epl+|\WE^{1-\tau}|\emi$ when $E=F_1$.
\endproclaim
\demo{Proof} Note that $1-\tau$ induces an isomorphism from
$\WE/\WF$ to $\WE^{1-\tau}$. Thus from Cor. 3.3, we have $E=F_1$ if
and only if $2$ divides $|\WE/\WF|$, which holds if and only if $2$
divides $|\WE^{1-\tau}|$. Also note that the greatest common divisor
of the orders of the subgroups
 $\WE^{1+\tau}$ and $\WE^{1-\tau}$ of the cyclic group $\WE$
 is the order of their intersection. This intersection clearly has exponent 2
and is cyclic. From Lemma 3.1 and Proposition 3.4, we know that the first group
$\WE^{1+\tau}$ has order divisible by 4.  We conclude that
this greatest common divisor is 2 if $E=F_1$ and 1 otherwise.

The elements $\alpha^+=|\WE^{1+\tau}|(1+\tau)$
and $\alpha^-=|\WE^{1-\tau}|(1-\tau)$ clearly annihilate $\WE$. When
 $E\neq F_1$, we show that these generate
$|\WE|$ and $\tau-N$ for some $N$, which in turn generate
$\Ann_\R(\WE)$. In this case, the integers $|\WE^{1+\tau}|$ and
$|\WE^{1-\tau}|$ are relatively prime, so an integer linear
combination of $\alpha^+$ and $\alpha^-$  yields the desired element
of the form $\tau-N$. Also $|\WE|/2$ is a common multiple of
$|\WE^{1+\tau}|$ and $|\WE^{1-\tau}|$, as
$|\WE|/\|WE^{1+\tau}|=|\WE^-|$, where $-1\in \WE^-$; while
$|\WE|/|\WE^{1-\tau}|=|\WE^G|=|\WF|=8m$ by Lemma 3.1 again. So
$\alpha^+$ and $\alpha^-$ generate $(|\WE|/2)(1+\tau + 1 -\tau)=|\WE|$.

Now suppose that $E=F_1$. We first show that
$\gamma=|\WE^{1+\tau}|\epl+|\WE^{1-\tau}|\emi=a\epl + b\emi \in \Ann_\R(\WE)$.
We have seen that the greatest common divisor of the
two orders $a$ and $b$ here is 2. If $\omega=\omega_1^2\in\WE^2$, then
$\omega^\gamma=(\omega_1^{1+\tau})^{|\WE^{1+\tau}|}(\omega_1^{1-\tau})^{|\WE^{1-\tau}|}=1$.
On the other hand, if $\omega \notin \WE^2$, then
its image generates $\WE/\WE^2$, and thus the image of
$\omega^{1+\tau}$ generates
$\WE^{1+\tau}/(\WE^{1+\tau})^2$,
and finally $\omega^{|\WE^{1+\tau}|\epl}$ generates
$\WE^{(1+\tau)|\WE^{1+\tau}|/2}$ which is cyclic
of order 2 and therefore equals $\{-1\}$. Similarly
$\omega^{|\WE^{1-\tau}|\emi}=-1$. Thus $\omega^\gamma=-1\cdot -1=1$
Next we show that $\gamma$ generates both $\tau-N$ for some $N$, and $|\WE|$.
Since $\gcd(a/2,b/2)=1$ and
$b/2$ is odd while $a/2$ is even, we have $\gcd(a,\frac a2- \frac b2)=1$. This shows that $\tau-N$ for some $N$
is an integer linear combination of $(1+\tau)\gamma=a\tau+a$ and
$\gamma = (\frac a2-\frac b2)\tau +(\frac a2+ \frac b2)$.
Also, $\frac{|\WE|}2$ is a common multiple of $a$ and $b$, as seen above, and
hence $|\WE|=\frac{|\WE|}{2}(1+\tau + 1-\tau)$
is generated by $\gamma(1+\tau)=a(1+\tau)$ and $\gamma(1-\tau)=b(1-\tau)$.
\qed
\enddemo

\heading{V. Properties of Fitting Ideals}
\endheading

In the following Proposition, we list the standard properties of
Fitting ideals which we will need. So suppose that $A$ is a
commutative Noetherian ring and $M$ is a finitely generated
$A$-module. We denote the Fitting ideal of $M$ over $A$ as
$\Fit_A(M)$. It is the ideal of $A$ generated by the determinants
of all square matrices representing relations among a set of
generators of $M$.

\proclaim{Proposition 5.1}
\item{1.} If $M$ is a cyclic $A$-module, then $\Fit_A(M)=\Ann_A(M)$.
\item{2.} Given a morphism of commutative Noetherian rings $f:A\rightarrow B$,
we have $\Fit_B(M\otimes_A B)=f(\Fit_A(M))B$. In particular,
\itemitem{a.} If $A\subset B$, then $\Fit_B(M\otimes_A B)=\Fit_A(M)B$.
\itemitem{b.} If $I$ is an ideal of $A$, then $\Fit_{A/I}(M/IM)=(\Fit_A(M)+I)/I$.
\itemitem{c.} If $T$ is a multiplicatively closed set in $A$ and $A_T=T^{-1}A$
(resp. $T^{-1}M=M_T$) is the corresponding
ring (resp. module) of fractions, then  $\Fit_{A_T}(M_T)=(\Fit_A(M))A_T$.
\item{3.} If $A=\Z$ and $M$ is finite, then $\Fit_{\Z}(M)=|M|\Z$.
\item{4.} If $A=A_1\oplus A_2$ and correspondingly $M=M_1 \oplus M_2$,
then $\Fit_A(M)=\Fit_{A_1}(M_1)\oplus \Fit_{A_2}(M_2)$
\endproclaim
\demo{Proof}
Proofs of parts 1--3 may be found in [\Eis]. The proof of part 4, and indeed
every part, is a straightforward exercise.
\qed
\enddemo

 We will also use a more specialized result on Fitting ideals,
variations of which are found in [\Gr] and [\Sch].

\proclaim{Proposition 5.2} Suppose that $G$ is a finite abelian
group and that $M$ is a finite $\Z[G]$-module. If $M$ is
cohomologically trivial, then $\Fit_{\Z_{(p)}[G]}(M_{(p)})$ is
principal, for any prime $p$.
\endproclaim
\demo{Proof} It is a standard fact [\AW, Thm. 9] that
$M$ is cohomologically trivial if and only if its projective dimension
as a $\Z[G]$-module is less than or equal to 1.
This means that we have a resolution of $M$ by projective modules
$$0\rightarrow Q \rightarrow P \rightarrow M\rightarrow 0,$$ and we may choose
$P$ to be free of some rank $n$.
Localization preserves exactness, so we get
$$0\rightarrow Q_{(p)} \rightarrow P_{(p)} \rightarrow M_{(p)}\rightarrow 0.$$
Now $\Z_{(p)}[G]$ is a semilocal ring (its maximal ideals correspond to those
of the finite ring $\Z/p\Z[G]$),  while $P_{(p)}$ is a free $\Z_{(p)}[G]$-module of rank $n$
and $Q_{(p)}$ is a projective submodule. We claim that it has constant rank $n$.
This is because if we localize the last exact sequence at a
minimal prime (necessarily contained in the zero-divisors) of $\Z_{(p)}[G]$,
the order of $M$ will be inverted and the middle term stays free of
rank $n$, while the term to the right becomes $0$ and hence the term to
the left must be free of rank $n$. Hence if we localize the projective module $Q_{(p)}$
at any prime of $\Z_{(p)}[G]$, it must be free, and further localization at a minimal
prime shows that the rank must be $n$. But a finitely generated projective
module of constant rank $n$ over a semilocal ring is free of rank $n$
(see eg. [\Eis, Exercise 4.13]). In this situation, the Fitting ideal
$\Fit_{\Z_{(p)}[G]}(M_{(p)})$
is clearly generated by the determinant of the single $n$ by $n$ matrix representing
the map in our exact sequence between the free module $Q_{(p)}$ and the free module
$P_{(p)}$.
\qed
\enddemo

\heading{VI. Fitting ideals in the maximal order of $\Q[G]$}
\endheading

We will have occasion to extend ideals from $\R=\Z[G]=\Z\oplus \Z \tau$ to
$\S=\Z \epl \oplus \Z \emi $, the maximal order in
$\Q[G]$. Note that $\S \cong \R/(1-\tau) \oplus \R/(1+\tau)$
under the obvious homomorphism, the inverse homomorphism
being $(\overline{r_1},\overline{r_2}) \rightarrow r_1\epl + r_2\emi$.
Also, $\Z\cong \R/(1\pm \tau)$
via the homomorphism sending each element of $\Z$ to its coset.
We use these isomorphisms to make identifications in the next lemma.

\proclaim{Lemma 6.1} Suppose that $M$ is a finite $\R$-module. Then the $\S$-ideal
$\Fit_{\R}(M)\S$ is generated by $|M^G|\epl +|M^-|\emi$.
\endproclaim
\demo{Proof} Note that $|M^G|=|M/M^{1-\tau}|$ since $M/M^G\cong M^{1-\tau}$,
and similarly $|M^-|=|M/M^{1+\tau}|$.
Using part 2a of Proposition 5.1, then part 4, and then part 3, we get
$$
\multline
\Fit_\R(M)\S=\Fit_\S(M\otimes_\R \S)=\Fit_{\R/(1-\tau)\oplus \R/(1+\tau)}(M/M^{1-\tau} \oplus M/M^{1+\tau})
\\=\Fit_{\R/(1-\tau)}(M/M^{1-\tau})\oplus \Fit_{\R/(1+\tau)}(M/M^{1+\tau})
\\
=(|M/M^{1-\tau}|\epl + |M/M^{1+\tau}|\emi)\S
=(|M^G|\epl + |M^-|\emi)\S
\qed
\endmultline
$$
\enddemo

Let $k^+=|\KE^G|$ (which equals $|\KF|$ by Proposition 2.2),
 and $k^-=|\KE^-|$.
Similarly let $w^+=|\WE^G|=|\WF|$ and $w^-=|\WE^-|$.

\proclaim{Proposition 6.2}
$$\FI\S=\(k^+\epl + k^- \emi\)\S$$ and
$$\Ann_{\R}(\WE)\S=\Fit_\R(\WE)\S=\(w^+\epl + w^- \emi \)\S$$
\endproclaim
\demo{Proof}
Apply the lemma to $\KE$ and $\WE$, and use part 1 of Prop. 5.1.
\qed
\enddemo

Our goal of course is to compute $\FI=\Fit_\R(\KE)$; we will also
need $\Ann_\R(\WE)=\Fit_\R(\WE)$. Prop. 6.2 will
help us identify the former. The latter is known and easily done
directly. We will derive it in a form most useful to us in the next
section. That, combined with Proposition 6.2 and the evaluation of
$\theta_{E/F}^S(-1)$ will be enough to provide the odd part of our
main result (Theorem 1.3). To obtain the 2-part of Theorem 1.3, we
will explicitly compute a key element of $\FI$ in Section IX.

\heading{VII. The Stickelberger Ideal
$\SI=\Ann_\R(\WE)\theta_{E/F}^S(-1)$}
\endheading

 Associated with each character $\chi$ of $G$, we have the Artin $L$-function with
Euler factors for primes in $S$ removed:

$${
L_{E/F}^S(s,\chi)=
\sum\Sb {\frak a} \ \text{integral} \\ ({\frak a},S)=1 \endSb
\dfrac{\chi(\sigma_{\fa})}{ {\N} {\frak a}^s}=
\prod\Sb \text{prime } \frak p \notin S \endSb
\bigl(1-\dfrac{\chi(\sigma_{\frak p})}{ {\N} {\frak p}^{s}}} \bigr)^{-1}.$$

These are related to the equivariant $L$-function by
$$\theta_{E/F}^S(s)=\sum_\chi L_{E/F}^S(s,\overline{\chi}) e_{\chi},
$$
where $\overline{\chi}$ is the complex conjugate of $\chi$,
as the two sides of this equation agree after multiplication by each idempotent
$e_\chi=\frac{1}{|G|}\sum_{\sigma\in G}\chi(\sigma)\sigma^{-1}$.
In our case we just have the trivial character $\chi_0$ and the
non-trivial character $\chi_1$, giving
$$
\multline \theta_{E/F}^S(s)=L_{E/F}^S(s,\chi_0)\epl
+L_{E/F}^S(s,\chi_1)\emi \\
 =\zeta_F^S(s)\epl +
\dfrac{\zeta_E^S(s)}{\zeta_F^S(s)}\emi,
\endmultline
$$
by standard properties of Artin $L$-functions. Thus
$$
\theta_{E/F}^S(-1)=\zeta_F^S(-1)\epl + \dfrac{\zeta_E^S(-1)}{\zeta_F^S(-1)}\emi.
$$  Referring now to  Conjecture 1.1 and Theorem 1.2, we clearly obtain the following.

\proclaim{Proposition 7.1} Assuming the Birch-Tate conjecture for
$E$ and for $F$, we have
$$(-1)^{|S|}\theta_{E/F}^S(-1)=\dfrac{|\KF|}{w_2(F)}\epl + (-1)^{|S_E|}\dfrac{w_2(F)}{w_2(E)}\dfrac{|\KE|}{|\KF|}\emi$$
Equality of the odd parts holds unconditionally.
\endproclaim

\proclaim{Lemma 7.2} Assuming the Birch-Tate conjecture for $E$ and
$F$, we have $w^+\epl\theta_{E/F}^S(-1)=\pm k^+\epl$ and $ w^- \emi
\theta_{E/F}^S(-1)=\pm k^- \emi$. Equality of odd parts holds
unconditionally.
\endproclaim
\demo{Proof} Using Proposition 7.1, we see that
$$
\multline \pm  w^+\epl \theta_{E/F}^S(-1)
=
\(w^+\epl\)\(\dfrac{k^+}{w^+}\epl \pm
\dfrac{w^+}{w_2(E)}\dfrac{|\KE|}{k^+}\emi\)
\\
 = k^+ \epl.
\endmultline
$$
Similarly,
$$
\multline \pm w^-\emi\theta_{E/F}^S(-1) =
\(w^-\emi\)\(\dfrac{k^+}{w^+}\epl \pm
\dfrac{w^+}{w_2(E)}\dfrac{|\KE|}{k^+}\emi\)
\\
 =
\pm \dfrac{w^+w^-}{w_2(E)}\dfrac{|\KE|}{k^+} \emi .
\endmultline
$$
Here $w_2(E)/w^-=|\WE^{1+\tau}|$, so
$$\dfrac{w^+w^-}{w_2(E)}=|\WE^G|/|\WE^{1+\tau}|
=|\KE^G|/|\KE^{1+\tau}|,$$ by Prop. 3.4 and Prop. 2.3. By Prop 2.2,
we also have $|\KE^G|=|\KF|=k^+$. Making these substitutions above
leads to the equality
$$
\pm (w^-\emi)\theta_{E/F}^S(-1) = \(|\KE|/|\KE^{1+\tau}|\) \emi= k^-
\emi,
$$
as desired. This proof works for the odd parts without assuming the
Birch-Tate conjecture.
\qed
\enddemo

 \proclaim{Proposition 7.3} Assuming the Birch-Tate
conjecture for $E$ and $F$, we have $\SI \S = \FI \S$. The equality
$\SI \S[1/2] = \FI \S[1/2]$ holds unconditionally.
\endproclaim
\demo{Proof} From Prop. 6.2 we have that  $\FI \S$
is generated by $k^+\epl + k^- \emi$ as an ideal in
$\S$. Similarly $\SI \S=\Ann_\R(\WE)\theta_{E/F}^S(-1)$ is
generated by $(w^+\epl + w^- \emi)\theta_{E/F}^S(-1)$. Thus it suffices to show that
 $k^+\epl + k^- \emi$ and $(w^+\epl + w^- \emi)\theta_{E/F}^S(-1)$
are associates in $\S$. This is easliy accomplished upon
multiplying by one of the units $\pm 1$ or
$\pm \tau$, and using Proposition 7.2. If we
first extend ideals to $\S[1/2]$, the same proof works
unconditionally.
\qed
\enddemo

\proclaim{Corollary 7.4} $\FI \R[1/2]=\SI \R[1/2]$. Consequently
$\FI_{(p)}=\SI_{(p)}$ for all odd primes $p$.
\endproclaim
\demo{Proof} Noting that $\R[1/2]=\S[1/2]$, we see that the first
equality is just a restatement of the unconditional part of
Proposition 7.3. The second equality follows upon extending ideals
from $\R[1/2]=\Z[1/2][G]$ to $\R_{(p)}=\Z_{(p)}[G]$.
\qed
\enddemo

 Cor. 7.4 establishes the unconditional part of Theorem 1.3.
Our goal now is to show that $\FI_{(2)}=\SI_{(2)}$, assuming
the 2-part of the Birch-Tate conjecture.

\heading{VIII. Ideals in $\R_{(2)}=\Z_{(2)}[G]$}
\endheading

 Since $G$ is a 2-group, $\R_{(2)}=\Z_{(2)}[G]$ is a local ring. This is easily seen as
it is an integral extension of the discrete valuation ring $\Z_{(2)}$, so any maximal ideal must contain
2, and hence must also contain $1-\tau$, whose square $2(1-\tau)$ is in the ideal generated by 2. The ideal $\m=(2,1-\tau)=(1+\tau,1-\tau)$ is of index 2, hence is the unique maximal ideal.

 We will need to consider the relationship between $\R_{(2)}$ and the overring $\S_{(2)}=\Z_{(2)}\epl\oplus \Z_{(2)}\emi\cong
\Z_{(2)}\oplus \Z_{(2)}$, which is clearly a principal ideal ring.

\proclaim{Lemma 8.1} The group of units $\S_{(2)}^\times $ of $\S_{(2)}$ equals the group of units of $\R_{(2)}$.
\endproclaim
\demo{Proof}
$$
\multline
\S_{(2)}^\times=\Z_{(2)}^\times\epl \oplus \Z_{(2)}^\times \emi=
(1+2\Z_{(2)})\epl \oplus (1+2\Z_{(2)})\emi=
\\
(\epl+\emi)+\Z_{(2)}(1+\tau)+\Z_{(2)}(1-\tau)
\\
=1+\m =\R_{(2)}^\times
\qed
\endmultline
$$
 \enddemo

\proclaim{Lemma 8.2} Suppose that $I$ is an ideal of finite index in
$\R_{(2)}$. Let $\alpha$ be a generator in the principal ideal ring $\S_{(2)}$
for the extended ideal $I\S_{(2)}$. Then $\alpha \in I$.
\endproclaim
\demo{Proof}
  Modifying by units, we may assume that $\alpha=2^i\epl + 2^j \emi$.
Since $\S=\R\epl\oplus \R\emi$, we have that
$2^i\epl + 2^j\emi$ generates $I\S=I\epl\oplus I\emi$.
Thus $2^i\epl$ generates $I\epl$ over $\Z_{(2)}$ and
$2^j\emi$ generates $I\emi$.
Consequently
$2^i\epl +2^{j'}\emi \in I$ for some $j'\ge j$
and $2^{i'}\epl + 2^j \emi \in I$ for some $i'\ge i$.
Multiplying the second element by $(1-\tau)\in \R$,
we also find that $2^{j+1}\emi \in I$.
If $j'=j$, then we have $2^i\epl +2^{j}\emi \in I$
as desired. If not, then $2^i\epl +2^{j'}\emi \in I$
and $2^{j+1}\emi \in I$ generate $2^i\epl \in I$.
Combined with $2^{i'}\epl + 2^j \emi \in I$, this generates
$2^j\emi \in I$, and finally $2^i\epl + 2^j \emi \in I$.
\qed
\enddemo

{\bf Remark 8.3} In the setting of Lemma 8.2, one can in fact show
that $I$ is either the principal ideal generated by $\alpha$ or the
non-prinicipal ideal generated by $\alpha \epl$ and $\alpha \emi$.

\heading{IX. Computation of $\FI_{(2)}$}
\endheading

  Assuming the 2-part of the Birch-Tate conjecture for $E$ and $F$, we
  now perform the computation of $\FI_{(2)}$, separating it into two different cases.

  \proclaim{Proposition 9.1} Suppose that $E=F_1$, and assume that
   the 2-part of the Birch-Tate conjecture holds for $F$ and for  $E$.
Then $\FI_{(2)}=\SI_{(2)}$.
\endproclaim
\demo{Proof}
 Since  $E=F_1$, Proposition 2.3 and Proposition 5.2 show that
 $\FI_{(2)}$ is a principal ideal in $\R_{(2)}$.
 Let $\alpha\in \R_{(2)}$ be a generator for $\FI_{(2)}$. Similarly,
 Proposition 3.4 and Proposition 5.2 show that $\Ann_R(\WE)\R_{(2)}$ and hence
 $\SI_{(2)}=\theta_{E/F}^S(-1)\Ann_R(\WE)\R_{(2)}$ is principal.
 We can also see this more explicitly from Prop. 4.1.
Let $\beta\in \R_{(2)}$ be a generator for $\SI_{(2)}$. From
Proposition 7.3, we deduce that $\FI_{(2)}$ and $\SI_{(2)}$ extend
to the same ideal in $\S_{(2)}$. Thus $\alpha$ and $\beta$ are
associates in $\S_{(2)}$. By Lemma 8.1, they are associates in
$\R_{(2)}$, and thus $\FI_{(2)}=\SI_{(2)}$.
\qed
\enddemo

 We turn now to the case of $E\neq F_1$, in which
$|\KE^-/\KE^{1-\tau}|=2$, by Prop. 2.3. We will use a lemma
to simplify the computation.

\proclaim{Lemma 9.2} Suppose that $M$ is a finite abelian 2-group, and
$m\in M$ is an element of order 2. Then $M$ is isomorphic to a direct
product of non-trivial cyclic groups, one of which contains $m$.
\endproclaim
\demo{Proof} We may assume that $M=\oplus_{i=1}^r \Z/2^{c_i}\Z$,with
$1\le c_i\le c_{i+1}$ for each $i$ between 1 and  $r-1$, inclusive.
We denote the standard generators by
 $e_1=(\overline{1},\overline{0},\dots,\overline{0}),
\cdots,e_r=(\overline{0},\dots,\overline{0},\overline{1})$.
Then necessarily $m=(2^{c_1-1}\ep_1,2^{c_2-1}\ep_2,\dots)$,
with $\ep_i=\overline{1}$ or $\overline{0}$, and not all $\ep_i=\overline{0}$.
Let $t$ be the smallest index for which $\ep_t\neq \overline{0}$.
Put $e_t'=(\overline{0},\cdots,\overline{0},\ep_t=\overline{1},
2^{c_{t+1}-c_t}\ep_{t+1},\dots,2^{c_r-c_t}\ep_{r})$, so that
$2^{c_t-1} e_t'=m$. Then we can see that
$M$ is the direct product of the subgroups generated by the elements
\newline
$e_1,\dots,e_{t-1},e_t',e_{t+1},\dots,e_r$; and this gives the desired
decomposition of $M$.
\qed
\enddemo

\proclaim{Proposition 9.3} If $E\neq F_1$, then
$k^+\epl\in \FI_{(2)}$.
\endproclaim
\demo{Proof} First, $\FI_{(2)}=\Fit_{\R_{(2)}}(\KE_{(2)})$, by part
2c) of Prop. 5.1. We know by Proposition
2.3 that $\KE^-/\KE^{1-\tau}$ has order 2, and it follows by
decomposing $\KE$ into a 2-part and an odd part that
 $\KE_{(2)}^-/\KE_{(2)}^{1-\tau}$ has order 2.
Now apply Lemma 9.2 with
\newline
 $M$=$\KE_{(2)}/\KE_{(2)}^{1-\tau}$ and $m$ as
a generator of the subgroup
\newline
$\KE_{(2)}^-/\KE_{(2)}^{1-\tau}$ of order
2. Let $\overline{\gamma_i}$, $1\le i\le r$ be the generators of
$\KE_{(2)}/\KE_{(2)}^{1-\tau}$ obtained from the lemma, with each
$\gamma_i \in \KE_{(2)}$. Then $\overline{\gamma_i}$ has order
$d_i=2^{c_i}$, a positive power of 2 for each $i$, and $\prod_i d_i=
|\KE_{(2)}/\KE_{(2)}^{1-\tau}|$. Denote this order by $k^+_{(2)}$.
Again the decomposition of $\KE$ into an odd part and a 2-part shows
that $k^+_{(2)}$ is the 2-part of $|\KE/\KE^{1-\tau}|$, which equals
$|\KF|=|\KE^G|=k^+$ by Thm. 2.1 and Prop. 2.2. Thus $k^+$ and
$k^+_{(2)}$ are associates in $\Z_{(2)}\subset \R_{(2)}$ and it
suffices to show that $k^+_{(2)}\epl\in \Fit_{\R_{(2)}}(\KE_{(2)})$.

The elements $\overline{\gamma_i}$, for $i=1,\dots, r$ generate
$\KE_{(2)}/\KE_{(2)}^{1-\tau}$ as a group, hence as an
$\R_{(2)}$-module. But $1-\tau \in \m$, the maximal ideal of the
local ring $\R_{(2)}$, and hence by Nakayama's lemma, the elements
$\gamma_i$ for $i=1,\dots, r$ generate $\KE_{(2)}$ as an
$\R_{(2)}$-module. Finding $r$ relations among these $r$ generators
will provide an element of the Fitting ideal in question by taking
the determinant of the corresponding $r$-by-$r$ relations matrix.

The fact that the elements $\gamma_i$, for $i=1,\dots, r$ generate
$\KE_{(2)}$ over $\R_{(2)}$ implies that the elements
$\gamma_i^{1-\tau}$ for $i=1,\dots, r$ generate
$\KE_{(2)}^{1-\tau}$. For each $i$, $\overline{\gamma_i}$ has order
$d_i$, and this means that $\gamma_i^{d_i}\in \KE_{(2)}^{1-\tau}$.
Multiplying by $\gamma_i^{(\tau-1)d_i/2}\in \KE_{(2)}^{1-\tau}$ gives
us $\gamma_i^{(1+\tau)d_i/2}\in \KE_{(2)}^{1-\tau}$ for each $i$. As
the elements $\gamma_i^{1-\tau}$ for $i=1,\dots,r$ generate this
last $\R_{(2)}$-module, we have that, for each $i$,
$$\gamma_i^{(1+\tau)d_i/2}=\prod_{j=1}^r(\gamma_j^{1-\tau})^{b_{ij}},$$
for some elements $b_{ij}\in\R_{(2)}$. Let $D$ be the $r$ by $r$
diagonal matrix of the $d_i$, and $B$ be the matrix of the $b_{ij}$.
Then  the last equation shows that $\epl D- (1-\tau)B$ is a
relations matrix for the generators $\gamma_i$ of $\KE_{(2)}$. Hence
$$\delta=\det\bigl(\epl D- (1-\tau)B\bigr)\in
\Fit_{\R_{(2)}}(\KE_{(2)})=\FI_{(2)}.$$

Modulo $(1-\tau)$, we have $\delta \equiv \det(\epl D)=(\epl)^r
\det(D)=\epl \prod_i d_i =  \epl k^+_{(2)}$. Before considering
$\delta$ modulo $1+\tau$, recall that the application of Lemma 9.2
also gives us that $\overline{\gamma_t}^{d_t/2}$ generates
$\KE_{(2)}^-/\KE_{(2)}^{1-\tau}$. It follows that $\gamma_t^{d_t/2}\in
\KE_{(2)}^-$, and consequently $\gamma_t^{(d_t/2)(1+\tau)}$ is
trivial. Hence we may choose $b_{tj}=0$ for all $j$, making
$\det(B)=0$. Now we see that $\delta\equiv \det(-2B)=0 $ modulo
$(1+\tau)$. Thus
$\delta - k^+_{(2)}\epl \in
(1+\tau)\cap (1-\tau)=(0)$ in $\R_{(2)}$. We
conclude that $k^+_{(2)}\epl=\delta \in \FI_{(2)}$, and we know that
this implies the desired result.
\qed
\enddemo

\proclaim{Proposition 9.4} Suppose that $E\neq F_1$, and assume that
the Birch-Tate conjecture holds for $F$ and $E$. Then
$\FI_{(2)}=\SI_{(2)}$.
\endproclaim
\demo{Proof}
By Proposition 6.2, $k^+\epl+k^-\emi$ generates
$\FI \S$ as an ideal of $\S$, and hence this same element
generates the extended ideal $\FI \S_{(2)}$ as an ideal of
$\S_{(2)}$. Proposition 8.2 then implies that
$k^+\epl+k^-\emi \in \FI_{(2)}$. On the other hand,
Proposition 9.3 gives $k^+\epl \in \FI_{(2)}$.
Subtracting, we see that both $k^+\epl$ and $k^-\emi$
lie in $\FI_{(2)}$. Thus
$$\multline
\FI_{(2)} \supset
k^+\epl \R_{(2)} + k^-\emi \R_{(2)}
\\
=(k^+\epl + k^-\emi)(\R_{(2)}\epl +\R_{(2)}\emi)
\\
=(k^+\epl + k^-\emi)\S_{(2)}
=\FI_{(2)}\S_{(2)}\supset \FI_{(2)},
\endmultline
$$
by Prop. 6.2 again. We conclude that $\FI_{(2)}=k^+\epl \R_{(2)} +
k^-\emi \R_{(2)}$. Turning to $\SI_{(2)}$, we have from Proposition
4.1 that $\Ann_\R(\WE)$ is generated by $|\WE^{1+\tau}|(1+\tau)$ and
$|\WE^{1-\tau}|(1-\tau)$. By Proposition 2.3, these generators may
also be written as $w^+ \epl$ and $w^- \emi$. Hence
$\SI=\theta_{E/F}^S(-1)\Ann_\R(\WE)$ is generated by $w^+\epl
\theta_{E/F}^S(-1)$ and $w^- \emi \theta_{E/F}^S(-1)$. Using
Lemma 7.2, we see that these generators equal
$\pm k^+\epl$ and $\pm k^-\emi$. Consequently $\SI$ is generated by
$k^+\epl$ and $k^-\emi$, so
$$\SI_{(2)}=k^+\epl\R_{(2)}+k^-\emi\R_{(2)}=\FI_{(2)}.\qed$$
\enddemo

Now we can complete the proof of Theorem 1.3.
Cor. 7.4 proves the unconditional statement in Theorem 1.3,
and shows that $$\FI_{(p)}=\SI_{(p)}\subset \Z_{(p)}[G]\subset \Q[G],$$ for each odd prime $p$.
Assuming the Birch-Tate conjecture for $E$ and $F$, Prop. 9.1 and Prop. 9.4 show that
$\FI_{(2)}=\SI_{(2)}.$ Note that for any ideal $I$ of $\R$, the intersection over all rational primes $\cap_p I_{(p)}=I$, because
for each $\alpha \in \cap_p \FI_{(p)}$, the $\Z$-ideal $\{a\in\Z: a\alpha \in I\}$ is contained in no prime ideal of $\Z$, hence equals $\Z$. Thus
$$\FI=\cap_{p} \FI_{(p)}
=\cap_{p} \SI_{(p)}=\SI. \qed$$

\heading{ACKNOWLEDGEMENTS}
\endheading
We thank Cristian Popescu and David Solomon for helpful discussions on the subject matter of this paper.
We thank Jason Price for useful comments on the exposition.

\enddocument
\end